\documentclass[12pt]{article}
\usepackage{amsmath,amsfonts,amssymb,rotating,colordvi}
\usepackage{colordvi}
\def\qed{\nopagebreak\hfill{\rule{4pt}{7pt}}}

\newtheorem{theo}{Theorem}[section]

\newtheorem{conj}[theo]{Conjecture}

\numberwithin{equation}{section}

\newdimen\Squaresize \Squaresize=11pt
\newdimen\Thickness \Thickness=0.7pt
\def\Square#1{\hbox{\vrule width \Thickness
   \vbox to \Squaresize{\hrule height \Thickness\vss
    \hbox to \Squaresize{\hss#1\hss}
   \vss\hrule height\Thickness}
\unskip\vrule width \Thickness} \kern-\Thickness}

\def\Vsquare#1{\vbox{\Square{$#1$}}\kern-\Thickness}

\def\moins{\raise 1pt\hbox{{$\scriptstyle -$}}}
\begin{document}

\begin{center}
{ \large\bf Proof of Lassalle's Positivity Conjecture

 on Schur Functions}
\end{center}

\begin{center}
William Y. C. Chen$^{1}$, Anne X. Y. Ren$^{2}$, Arthur
L. B. Yang$^{3}$\\[6pt]

$^{1,2,3}$Center for Combinatorics, LPMC-TJKLC\\
Nankai University, Tianjin 300071, P. R. China

Email: $^{1}${\tt
chen@nankai.edu.cn}, $^{2}${\tt renxy@nankai.edu.cn}, $^{3}${\tt yang@nankai.edu.cn}
\end{center}

\noindent\textbf{Abstract.}
In the study of Zeilberger's conjecture on an integer sequence related to the Catalan numbers, Lassalle proposed the following conjecture. Let $(t)_n$ denote the rising factorial, and let $\Lambda_{\mathbb{R}}$ denote the algebra of symmetric functions with real coefficients. If $\varphi$ is the homomorphism from $\Lambda_{\mathbb{R}}$ to $\mathbb{R}$ defined by
$\varphi(h_n)={1}/{((t)_nn!)}$ for some $t>0$, then for any Schur function $s_{\lambda}$, the value $\varphi(s_{\lambda})$ is positive. In this paper, we provide an affirmative answer to Lassalle's conjecture  by using the Laguerre-P\'olya-Schur theory of multiplier sequences.

\noindent \emph{AMS Classification 2010:} Primary 05E05; Secondary 26C10.

\noindent \emph{Keywords:}  symmetric function, Schur function, multiplier sequence, totally positive sequence.

\section{Introduction}

The objective of this paper is to prove a positivity conjecture on Schur functions, which was proposed by Lassalle \cite{las2012} in the study of two combinatorial sequences related to the Catalan numbers.

Let us begin with an overview of Lassalle's conjecture. Let
$$C_n=\frac{1}{n+1}\binom{2n}{n}$$ denote the $n$-th Catalan number. Lassalle \cite{las2012} introduced a sequence of numbers $A_n$ for $n\geq 1$, which are recursively defined by
$$(-1)^{n-1}A_n=C_n+\sum_{j=1}^{n-1}(-1)^j \binom{2n-1}{2j-1}A_jC_{n-j},$$
with the initial value $A_1=1$. He proved that the sequence $\{A_n\}_{n\geq 2}$ is positive and increasing.
Josuat-Verg$\mathrm{\grave{e}}$s \cite{verges2012} found a combinatorial interpretation of $A_n$ in terms of connected matchings in the study of cumulants of the $q$-semicircular law.
Zeilberger further conjectured that the numbers $\{2A_n/C_n\}_{n\geq 2}$ also form an increasing sequence of positive integers. Lassalle \cite{las2012}  proved Zeilberger's conjecture. An alternative proof was given by Amdeberhan, Moll and Vignat \cite{amv2012} using a probabilistic approach.

By using the theory of symmetric functions, Lassalle \cite{las2012} gave a direct proof of the positivity and the monotonicity of $\{A_n\}_{n\geq 2}$, although these two properties can be deduced from Zeilberger's conjecture. For the notation and terminology on symmetric functions, see Stanley \cite{stanley1999}.
Lassalle's proof  involves the following specialization of symmetric functions.  Let $\mathbb{R}$ be the field of real numbers, and let $\Lambda_{\mathbb{R}}$ be the algebra of symmetric functions with real coefficients.
It is well known that the complete symmetric functions $h_n$ $(n\geq 0)$ are algebraically independent and  $\Lambda_{\mathbb{R}}$ is
generated by $h_n$. Thus any homomorphism $\varphi$ from $\Lambda_{\mathbb{R}}$ to $\mathbb{R}$ is uniquely determined by the values $\varphi(h_n)$. Lassalle's specialization is given by
\begin{align}\label{eq-specialization}
\varphi(h_n)=\frac{1}{((t)_nn!)},
\end{align}
where $t>0$ and $(t)_n=t(t+1)\cdots(t+n-1)$. Lassalle proved that this specialization satisfies
$$\varphi((-1)^{n-1}p_n)>0 \qquad \mbox{ and }\qquad \varphi(e_n)>0,$$
where $p_n$ and $e_n$
denote the $n$-th power sum and the $n$-th elementary symmetric function respectively.

Note that both $h_n$ and $e_{n}$ are special cases of the Schur functions. Based on the positivity of $\varphi(h_n)$ and $\varphi(e_n)$, Lassalle further considered the specialization of a general Schur function. Recall that an integer partition $\lambda$ is a weakly decreasing sequence $(\lambda_1,\lambda_2,\ldots,\lambda_\ell)$ of nonnegative integers. A skew partition is a pair of partitions $(\lambda,\mu)$ of the same length $\ell$ such that $\mu_i\leq \lambda_i$ for each $1\leq i\leq \ell$, denoted by $\lambda/\mu$. A skew Schur function $s_{\lambda/\mu}$, indexed by a skew partition $\lambda/\mu$, is given by
\begin{align}\label{eq-jt}
s_{\lambda/\mu}=\det(h_{\lambda_i-\mu_j-i+j})_{i,j=1}^{\ell},
\end{align}
where  $h_k$ is set to be zero if $k<0$. This formula is known as the Jacobi-Trudi identity for skew Schur functions.
If $\mu$ is the empty partition, then $s_{\lambda/\mu}$ is called a Schur function of shape $\lambda$, denoted by $s_{\lambda}$.
Lassalle posed the following conjecture.

\begin{conj}[{\cite[p.  933]{las2012}}]\label{conj-las}
Let $\varphi\colon \Lambda_{\mathbb{R}}\rightarrow \mathbb{R}$  be the specialization of $h_n$  given by \eqref{eq-specialization}.  Then $\varphi(s_{\lambda})$ is positive for any Schur function $s_{\lambda}$.
\end{conj}

In this paper, we  give an affirmative answer to Conjecture \ref{conj-las}.
Our proof relies on the theory of total positivity  and the theory of multiplier sequences.

\section{Preliminaries}

In this section, we give an overview of some fundamental results on the theory of total positivity and the theory of multiplier sequences.
A real sequence $\{a_n\}_{n\geq 0}$ is said to be a totally positive sequence if all the minors of the infinite Toeplitz matrix $(a_{j-i})_{i,j\geq 1}$ are nonnegative, where we set $a_n=0$ for $n<0$. The following  representation theorem  was conjectured by Schoenberg and proved by Edrei \cite{edr1953}, see also Karlin \cite{karlin1968}.

\begin{theo}[{\cite[p.  412]{karlin1968}}]\label{thm-total} Let $\{a_n\}_{n\geq 0}$ be a sequence of real numbers with $a_0=1$. Then $\{a_n\}_{n\geq 0}$ is totally positive if and only if its generating function
$$f(x)=\sum_{n\geq 0}a_n x^n$$
 has the form
\begin{align}\label{eq-total}
e^{\theta x}\frac{\prod_{i\geq 1}(1+\rho_ix)}{\prod_{i\geq 1}(1-\delta_ix)},
\end{align}
where $\theta\geq 0, \rho_i\geq 0, \delta_i\geq 0$ for $i\geq 1$ and $\sum_{i\geq 1}(\rho_i+\delta_i)<\infty$.
\end{theo}

Based on the above   theorem, Karlin gave a necessary and sufficient condition to determine the strict positivity of a minor of the
Toeplitz matrix  $(a_{j-i})_{i,j\geq 1}$.

\begin{theo}[{\cite[p. 428]{karlin1968}}] \label{thm-strict}
Suppose that $\{a_n\}_{n\geq 0}$ is a totally positive sequence.
Let $\theta, \delta_i, \rho_i$ be defined as in \eqref{eq-total}.
Let $K$ be the number of positive entries $\delta_i$ and let $L$ be the number of positive entries $\rho_i$, where  $K$ and $L$ are
allowed to be infinity.
Let $I=(i_1,i_2,\ldots,i_r)$ and $J=(j_1,j_2,\ldots,j_r)$ be two increasing sequences of positive numbers. Let $T(I,J)$ be the minor of $(a_{j-i})_{i,j\geq 1}$ with the row indices $i_1,i_2,\ldots,i_r$ and column indices $j_1,j_2,\ldots,j_r$.
Then the following assertions hold:
\begin{itemize}
\item[(i)] For $\theta>0$, the minor $T(I,J)$ is positive if and only if $i_k\leq j_k$ for $1\leq k\leq r$;

\item[(ii)] For $\theta=0$ and $K>0$, the minor $T(I,J)$ is positive if and only if $$j_{k-K}-L<i_k\leq j_k$$
     for $1\leq k\leq r$.

\item[(iii)]
     For $\theta=0$ and $K=0$, the minor $T(I,J)$ is positive if and only if $$j_{k}-L\leq i_k\leq j_k$$
     for $1\leq k\leq r$.
\end{itemize}
\end{theo}

As pointed out by Craven and Csordas \cite{cc1996}, Theorem \ref{thm-total} is closely related to
P\'olya and Schur's transcendental characterization of multiplier sequences.
A multiplier sequence is defined to be a sequence $\{\gamma_n\}_{n\geq 0}$ of real numbers such that, whenever the polynomial with real coefficients
$$\sum_{n=0}^m a_n x^n$$
has only real zeros, the polynomial
$$\sum_{n=0}^m \gamma_n a_n x^n$$
also  has only real zeros.
P\'olya and Schur obtained the following transcendental characterization of multiplier sequences consisting of nonnegative numbers,  see also Levin \cite{levin1964}.

\begin{theo}[{\cite[p.  346]{levin1964}}] \label{thm-ps} A sequence  $\{\gamma_n\}_{n\geq 0}$ of nonnegative numbers with $\gamma_0=1$ is a multiplier sequence if and only if
$$f(x)=\sum_{n\geq 0} \frac{\gamma_n}{n!}x^n$$ is of the form
\begin{align}\label{eq-total-special}
e^{\theta x}{\prod_{i\geq 1}(1+\rho_ix)},
\end{align}
where $\theta\geq 0, \rho_i\geq 0$ for $i\geq 1$ and $\sum_{i\geq 1} \rho_i<\infty$.
\end{theo}

P\'olya and Schur also gave an algebraic characterization of multiplier sequences.

\begin{theo}[{\cite[p. 345]{levin1964}}] \label{thm-ps-a} A sequence  $\{\gamma_n\}_{n\geq 0}$ of nonnegative numbers is a multiplier sequence if and only if for all $n\geq 0$ the Jensen polynomial
$$J_n(x)=\sum_{k=0}^n \binom{n}{k}\gamma_k x^k$$
has only nonpositive real zeros.
\end{theo}

As will be shown in Section 3,  Lassalle's conjecture is equivalent to the total positivity of the sequence $\left\{{1}/{((t)_{n}n!)}\right\}_{n\geq 0}$ for $t> 0$.
To prove the required total positivity, we
shall use a classic result of Laguerre on multiplier sequences. 
Suppose that $\{\gamma_n\}_{n\geq 0}$ is a nonnegative sequence with $\gamma_0=1$ such that its exponential generating function $f(x)$ is
 entire. By Theorems \ref{thm-total} and \ref{thm-ps}, the sequence $\{\gamma_n\}_{n\geq 0}$ is a multiplier sequence if and only if the sequence $\{\gamma_n/n!\}_{n\geq 0}$ is totally positive. Thus, the total positivity of the sequence $\left\{{1}/{((t)_{n}n!)}\right\}_{n\geq 0}$
  is a consequence of the following theorem, which can be traced back to Laguerre, see also Levin \cite{levin1964}.

\begin{theo}[{\cite[p. 341]{levin1964}}]\label{thm-lag}
For any $t>0$, the sequence $\left\{{1}/{(t)_{n}}\right\}_{n\geq 0}$ is a multiplier sequence.
\end{theo}

It is worth mentioning that the above theorem is closely related to the real-rootedness  of the Bessel functions.
In the study of zero-mapping transformations, Iserles, N\o rsett and Saff \cite{ins1991} showed that the exponential generating function $f(x)$ of the multiplier sequence $\left\{{1}/{(t)_{n}}\right\}_{n\geq 0}$ is given by
\begin{align}\label{eq-bessel}
f(x)=\sum_{n\geq 0} \frac{1}{(t)_{n}n!} x^n=(i\sqrt{x})^{(1-t)}\Gamma(t)J_{t-1}(2i\sqrt{x}),
\end{align}
where $\Gamma(t)$ is the Gamma function and $J_{t-1}(x)$ is the Bessel function of order $t-1$. From Hurwitz's theorem on the
 real-rootedness of the entire function $x^{1-t}J_{t-1}(x)$ \cite[p. 483]{watson1922}, we deduce
  that the sequence
   $\left\{{1}/{(t)_{n}}\right\}_{n\geq 0}$ satisfies the conditions of Theorem \ref{thm-ps}.

Theorem \ref{thm-lag} is also related to the real-rootedness  of the
Laguerre polynomials.
Iserles, N\o rsett and Saff \cite{ins1991} showed that the Jensen polynomials with respect to the multiplier sequence $\left\{{1}/{(t)_{n}}\right\}_{n\geq 0}$ are given by
\begin{align*}
J_n(x)=\sum_{k=0}^n \binom{n}{k}\frac{1}{(t)_k} x^k=\frac{n!}{(t)_n}L_n^{t-1}(-x),
\end{align*}
where $L_n^{t-1}(x)$ is the  Laguerre polynomial \cite{rain1967}. It is well known that $L_n^{t-1}(x)$ has only real zeros. Hence
Theorem \ref{thm-lag} can also be deduced from Theorem \ref{thm-ps-a}.

\section{Proof of Lassalle's conjecture} % Proof of Lassalle's conjecture

In this section, we prove the following theorem on the
positivity of the specialization of
skew Schur functions. It is clear that
Lassalle's conjecture is a consequence of Theorem 3.1.
However, it should be mentioned that Lassalle's conjecture
is in fact equivalent to Theorem 3.1 since each
skew Schur function is a linear combination of Schur
function with nonnegative coefficients. The reason for us to
use the formulation in terms of skew Schur functions lies in the
 connection between skew Schur functions and non-vanishing minors of the Toeplitz matrix $(a_{j-i})_{i,j\geq 1}$ corresponding to a totally positive sequence $\{a_n\}_{n\geq 0}$.

\begin{theo}\label{thm-main}
Suppose that $\varphi\colon \Lambda_{\mathbb{R}}\rightarrow \mathbb{R}$ is the specialization of $h_n$  given by \eqref{eq-specialization}. Then $\varphi(s_{\lambda/\mu})$ is positive for any skew Schur function $s_{\lambda/\mu}$.
\end{theo}

\noindent \emph{Proof of Theorem \ref{thm-main}.}
 We shall divide the proof into three steps. The first step is to transform Theorem \ref{thm-main} to a problem of determining the total positivity of certain sequence by using the Jacobi-Trudi identity for skew Schur functions. The second step is to prove the total positivity of this sequence by making a connection to Theorem \ref{thm-lag}. The third step is to prove the strict positivity of $\varphi(s_{\lambda})$ by using
Karlin's criterion for determining the strict positivity of a minor of $(a_{j-i})_{i,j\geq 1}$ for a totally positive sequence $\{a_n\}_{n\geq 0}$.

First, we prove that Theorem \ref{thm-main} is equivalent to the assertion that the sequence
$$\left\{\frac{1}{({(t)_{n}(n)!})}\right\}_{n\geq 0}$$
is totally positive. Let $T=(T_{i,j})_{i,j\geq 1}$ be the Toeplitz matrix corresponding to the above sequence, namely
$$
T_{i,j}=\left\{
\begin{array}{ll}
\frac{1}{(t)_{j-i}(j-i)!}, & \mbox{ if }i\leq j,\\[5pt]
0, & \mbox{ otherwise}.
\end{array}
\right.
$$
It suffices to show that every
$\varphi(s_{\lambda/\mu})$ occurs as a minor of $T$, and every minor of $T$
is equal to $\varphi(s_{\lambda/\mu})$ for some skew partition $\lambda/\mu$. It is worth mentioning that a similar situation
also happens to $G$-analogues of symmetric functions, as observed by  Stanley \cite[Corollary 2.2]{stanley1995}. The proof is a straightforward application of the Jacobi-Trudi identity.
Given a skew partition $\lambda/\mu$ with $\lambda=(\lambda_1,\lambda_2,\ldots,\lambda_{\ell})$ and $ \mu=(\mu_1,\mu_2,\ldots,\mu_{\ell}),$
let
\begin{align}
I &= (\mu_{\ell}+1,\mu_{\ell-1}+2,\ldots,\mu_1+\ell),\label{eq-i}\\
J &= (\lambda_{\ell}+1,\lambda_{\ell-1}+2,\ldots,\lambda_1+\ell).\label{eq-j}
\end{align}
Let $T(I,J)$ denote the minor of the  Toeplitz matrix $T$
 with row index set $I$ and    column index set $J$. By the Jacobi-Trudi identity \eqref{eq-jt}, we obtain
$$\varphi(s_{\lambda/\mu})=T(I,J).$$
In other words, the specialization
$\varphi(s_{\lambda/\mu})$ always occurs as a minor of $T$.
Conversely, any minor of $T$
is equal to $\varphi(s_{\lambda/\mu})$ for certain skew Schur function.  This completes the proof of the assertion of the first step.

We now proceed to prove the total positivity of the sequence $\{{1}/({(t)_{n}n!})\}_{n\geq 0}$.
Combing Theorems \ref{thm-ps} and \ref{thm-lag}, the generating
function
$$f(x)=\sum_{n\geq 0} \frac{1}{(t)_{n}n!}x^n$$
is entire and has the form \eqref{eq-total-special}. By Theorem \ref{thm-total}, we are led to total positivity.
Alternatively, we can use \eqref{eq-bessel} to prove the total positivity based on the real-rootedness of the Bessel functions.

Finally, we prove the strict positivity of the specialization $\varphi(s_{\lambda/\mu})$  for any skew partition $\lambda/\mu$.
By the total positivity of the sequence $\{{1}/({(t)_{n}n!})\}_{n\geq 0}$, we see that $\varphi(s_{\lambda/\mu})$ is nonnegative.
Suppose that $\varphi(s_{\lambda/\mu})=T(I,J)$ with $I$ and $J$ given by \eqref{eq-i} and \eqref{eq-j}.
By using Karlin's criterion, namely Theorem \ref{thm-strict}, we   are going to prove  the strict positivity of $T(I,J)$.
For the sequence $\{{1}/({(t)_{n}n!})\}_{n\geq 0}$, we need to consider the values of the parameters $K,L$ and $\theta$ which appear in Theorem \ref{thm-strict}. Since the generating function $f(x)$ is of the form \eqref{eq-total-special}, we see that $K=0$. If $\theta>0$, then we have $T(I,J)>0$,  since, for $1\leq k \leq \ell$,
$$i_k=\mu_{\ell+1-k}+k\leq \lambda_{\ell+1-k}+k=j_k.$$
If $\theta=0$, then we have $L=\infty$, since $f(x)$ is not a polynomial. By (iii) of Theorem \ref{thm-strict}, we have $T(I,J)>0$,  since the condition
$$j_{k-K}-L\leq i_k\leq j_k$$ is  satisfied for $1\leq k \leq \ell$. In both cases, we obtain $T(I,J)>0$, and hence $\varphi(s_{\lambda/\mu})>0$.
This completes the proof of the theorem. \qed

\vskip 3mm
\noindent {\bf Acknowledgments.} This work was supported by the 973 Project, the PCSIRT Project
of the Ministry of Education, and the National Science Foundation of China.

\end{document}